\documentclass{article}

\usepackage[english]{babel}
\usepackage{amsthm}
\usepackage{amsmath}
\usepackage{amssymb}

\usepackage{graphicx}
\usepackage{enumitem}
\usepackage[colorlinks]{hyperref}
\usepackage{caption}
\usepackage{subcaption}
\usepackage{pict2e,color,colordvi,amscd}
\usepackage{indentfirst}
\usepackage{tikz}
\usepackage{float}

\newtheorem{definition}{Definition}[section]
\newtheorem{theorem}[definition]{Theorem}
\newtheorem{lemma}[definition]{Lemma}

\newtheorem{corollary}[definition]{Corollary}

\newtheorem{conjecture}{Conjecture}

\theoremstyle{definition}

\theoremstyle{plain}

\title{Claw-free minimal bricks}
\author{
Guantao Chen\thanks{Department of Mathematics and Statistics, Georgia State University, Atlanta, GA 30303, USA. Email address: gchen@gsu.edu.}
\and
Yuying Ma\thanks{Department of Mathematics and Statistics, Georgia State University, Atlanta, GA 30303, USA. Email address: yma29@gsu.edu.}
}
\date{\relax}
\begin{document}

\maketitle

\begin{abstract}
A graph $G$ is called a \emph{brick} if it is $3$-connected and $G-\{u,v\}$ has a perfect matching for every pair of
distinct vertices $u,v\in V(G)$. A brick is \emph{minimal} if deleting any edge results in a graph that is no longer a brick. A graph is \emph{claw-free} if it contains no induced $K_{1,3}$. Plummer proved that every $3$-connected claw-free graph of even order is a brick.

In this paper, we completely characterize all claw-free minimal bricks. Apart from $K_4$ and the six-vertex wheel $W_5$, they are precisely the line graphs of graphs obtained from loopless $3$-edge-connected cubic multigraphs by subdividing every edge outside a matching of even size exactly once.

Consequently, every claw-free minimal brick of order $n$ contains at least $4n/5$ cubic vertices and at least $4n/5$ edges joining cubic vertices. Both bounds are sharp. This establishes the conjectures of Norine and Thomas and of Lov\'asz for claw-free minimal bricks in sharp quantitative forms.

We also show that, for a claw-free minimal brick of order at least $8$, deleting an edge leaves a matching covered graph if and only if the edge lies in a triangle.

\vspace{0.3cm}
\end{abstract}

\section{Introduction}
Throughout this paper, all graphs are finite and simple unless stated otherwise. For a graph $G$, let $|G|:=|V(G)|$ denote its \emph{order}. A \emph{matching} in $G$ is a set of pairwise vertex-disjoint edges, and it is \emph{perfect} if it covers every vertex of $G$. A nontrivial connected graph is \emph{matching covered} if every edge belongs to a perfect matching. A graph $G$ is \emph{bicritical} if $G-\{u,v\}$ has a perfect matching for every pair of distinct vertices
$u,v\in V(G)$. A \emph{brick} is a $3$-connected bicritical graph. Bricks play a fundamental role in matching theory as the nonbipartite building blocks in the tight-cut decomposition of matching covered graphs~\cite{LOVASZ1987}.

A brick is \emph{minimal} if deleting any edge results in a graph that is no longer a brick. Since every brick contains a spanning minimal brick, understanding the structure of minimal bricks is of particular interest.
A graph is \emph{claw-free} if it contains no induced $K_{1,3}$. Plummer~\cite{Plummer1994} proved that every $3$-connected claw-free graph of even order is a brick. On the other hand, a matching covered graph $G$ is \emph{minimally matching covered} if $G-e$ is not matching covered for every $e\in E(G)$. Zhang et al.~\cite{ZHANG2025} characterized all claw-free minimally matching covered graphs. The distinction between these two notions of minimality is substantial: Lov\'asz~\cite{LOVASZ1987} proved that $K_4$ and the triangular prism are the only bricks that are minimally matching covered.

Thus, while claw-free minimally matching covered graphs have been characterized, the problem of characterizing claw-free minimal bricks remains open. The main purpose of this paper is to determine the structure of this class and investigate several consequences of this characterization.

The \emph{line graph} of a graph $H$, denoted by $L(H)$, has vertex set $E(H)$, with two vertices adjacent if and only if the corresponding edges of $H$ share an endvertex. For $n\ge 3$, let $W_n$ denote the graph obtained from an $n$-cycle by adding a vertex adjacent to every vertex of the cycle. Recall that subdividing an edge once means replacing it with a path of length two.

\begin{theorem}\label{thm:main1-char}
A graph is a claw-free minimal brick if and only if it is $K_4$, $W_5$, or the line graph of a graph obtained from a loopless, $3$-edge-connected cubic multigraph by subdividing each edge outside a matching of even size (possibly empty) exactly once.
\end{theorem}

The only nonsimple multigraph in this construction consists of two vertices joined by three parallel edges. In this case, the matching of even size must be empty, so every edge is subdivided once, and the resulting line graph is the triangular prism. Moreover, Theorem~\ref{thm:main1-char} implies that the possible orders of claw-free minimal bricks are precisely the positive even integers other than $2$, $8$, and $14$. 

There are two longstanding conjectures on the distribution of cubic vertices in minimal bricks.

\begin{conjecture}[Norine and Thomas~\cite{NORINE2006}]
\label{conj:RT}
There exists a constant $\alpha > 0$ such that every minimal brick $G$ has at least $\alpha|G|$ cubic vertices.
\end{conjecture}

\begin{conjecture}[Lov\'asz (see~\cite{Yuliu2009})]
\label{con:L}
Every minimal brick contains two adjacent cubic vertices.
\end{conjecture}

Several partial results toward these conjectures are known. De Carvalho, Lucchesi and Murty~\cite{DECARVALHO2006} proved that every minimal brick has a cubic vertex. Norine and Thomas~\cite{NORINE2006} improved this lower bound to three, and Lin, Zhang and Lu~\cite{LZL2014} subsequently improved it to four. Feng~\cite{Feng2024} verified Conjecture~\ref{conj:RT} for minimal bricks with average degree at least $4.6$, while He and Lu~\cite{HE2024} proved that every solid minimal brick $G$ has at least $2|G|/5$ cubic vertices. More recently, Xue and Lu~\cite{XL2025thin} showed that every claw-free minimal brick $G$ of order at least $8$ has at least $3|G|/16$ cubic vertices.
Progress on Conjecture~\ref{con:L} has been more limited. Chen, Feng and Lu~\cite{CFL2026} proved that every minimal brick contains two cubic vertices at distance at most $2$, and verified the conjecture for minimal bricks with average degree at least $4.5$.

As a consequence of our structural characterization, we obtain the following sharp bounds on the cubic vertices of claw-free minimal bricks.

\begin{theorem}\label{thm:counting}
Let $G$ be a claw-free minimal brick of order $n$, and let $V_3$ be the set of its cubic vertices. Then
\[
|V_3|\ge \frac{4n}{5}
\qquad\text{and}\qquad
|E(G[V_3])|\ge |V_3|.
\]
Both inequalities are sharp. Moreover, if $G\neq W_5$, then $G[V_3]$ admits a perfect matching.
\end{theorem}

Theorem~\ref{thm:counting} improves the lower bound $3|G|/16$ of Xue and Lu~\cite{XL2025thin} for the number of cubic vertices to the sharp bound $4|G|/5$. Moreover, it shows that every claw-free minimal brick $G$ contains at least $4|G|/5$ edges joining cubic vertices. If $G \neq W_5$, at least $2|G|/5$ of these edges are pairwise vertex-disjoint. Thus both conjectures hold for claw-free minimal bricks in stronger quantitative forms.

We next consider three important classes of edges in matching theory: removable, $b$-invariant, and thin edges. An edge $e$ of a matching covered graph $G$ is \emph{removable} if $G-e$ is also matching covered. Note that every brick other than $K_4$ and the triangular prism has a removable edge.

To define the other two classes, recall that an edge cut is \emph{tight} if every perfect matching contains exactly one edge of the cut. By the characterization of Edmonds, Pulleyblank, and Lov\'asz~\cite{Edmonds1982}, bricks are precisely the nonbipartite matching covered graphs without nontrivial tight cuts. A \emph{brace} is a bipartite matching covered graph without nontrivial tight cuts. A \emph{tight-cut decomposition} of a matching covered graph is obtained by repeatedly contracting the shores of nontrivial tight cuts. This process yields a family of bricks and braces.
Lov\'asz~\cite{LOVASZ1987} proved that this family is unique up to multiple edges. In particular, the number of bricks in such a decomposition, denoted by $b(G)$, is well-defined. A removable edge $e$ of a brick $G$ is \emph{$b$-invariant} if $b(G-e) = b(G) = 1$.

The \emph{bicontraction} of a vertex $v$ of degree $2$ with distinct neighbors $u$ and $w$ consists of contracting $\{v, u, w\}$ into a single vertex. The \emph{retract} of a matching covered graph is obtained by repeatedly performing bicontractions until no vertex of degree $2$ with distinct neighbors remains. A removable edge $e$ of a brick $G$ is \emph{thin} if the retract of $G-e$ is a brick. Let $R(G)$, $I(G)$, and $T(G)$ denote the sets of removable, $b$-invariant, and thin edges of $G$, respectively. In general, $T(G)\subseteq I(G)\subseteq R(G)$.

Lov\'asz conjectured that every brick other than $K_4$, the triangular prism, and the Petersen graph
has a $b$-invariant edge. De Carvalho, Lucchesi and Murty~\cite{CLM2002I,CLM2002II} confirmed this conjecture and later proved that every brick other than these three graphs has a thin edge~\cite{DECARVALHO2006}. These existence results lead to the question of how many edges of each type a brick must contain.

For claw-free bricks $G$ of order at least $8$, Wu, Lu and Zhang~\cite{WLZ2024} showed that $|R(G)|\ge 5|G|/8$, while Xue and Lu~\cite{XL2025thin} showed that $|T(G)|\ge 3|G|/8$. For cubic claw-free bricks $G$ other than $K_4$ and the triangular prism, Zhang et al.~\cite{ZHANG2025} showed that $|R(G)|=|G|$.

For claw-free minimal bricks, our characterization yields a simple description of all three edge classes, together with sharp bounds on their sizes.

\begin{theorem}\label{thm:edge-classes}
Let $G$ be a claw-free minimal brick of order
at least $8$. Then
\[
R(G)=I(G)=T(G)
=\{e\in E(G): e\text{ lies in a triangle of }G\}.
\]
Moreover,
\(
|G|\le |R(G)|\le \frac{6|G|}{5},
\)
and both bounds are sharp.
\end{theorem}

The remainder of the paper is organized as follows. In the next section, we reduce the study of claw-free minimal bricks to that of claw-free minimally $3$-connected graph of even order. In Section~\ref{sec:pf-main}, we prove Theorem~\ref{thm:main1-char} and determine the possible orders of claw-free minimal bricks. In Section~\ref{sec:applications}, we apply the main theorem to show  Theorems~\ref{thm:counting} and~\ref{thm:edge-classes}.

We use the following notation throughout the paper. Let $G$ be a graph. Denote the maximum degree, minimum degree and independent number of $G$ by $\Delta(G)$, $\delta(G)$ and $\alpha(G)$, respectively. For $X\subseteq V(G)$ and $F\subseteq E(G)$, let $G-X$ denote the subgraph induced by $V(G)\setminus X$, and let $G-F$ denote the spanning subgraph of $G$ with edge set $E(G-F)=E(G)\setminus F$. When $X=\{x\}$ and $F=\{e\}$, we write $G-x$ for $G-\{x\}$, and $G-e$ for $G-\{e\}$. For terminology and notation not defined here, we follow Bondy and Murty~\cite{bondy2011}.

\section{Reduction to minimally $3$-connected}

We recall Tutte's $1$-factor theorem: a graph $G$ has a perfect matching if and only if $o(G-S)\le |S|$ for every $S\subseteq V(G)$, where $o(H)$ is the number of components of odd order in a graph $H$. Applying this theorem to $G-\{x,y\}$ for each pair of distinct vertices $x,y\in V(G)$ yields the following lemma. 

\begin{lemma}\label{lem:tutte-bicritical-fact}
A graph $G$ with $|G|\ge 2$ is bicritical if and only if
\(
    o(G-S)\le |S|-2
\)
for every $S\subseteq V(G)$ with $|S|\ge 2$.
\end{lemma}

We also recall two classical results of Plummer and Halin.

\begin{theorem}[Plummer~\cite{Plummer1994}]\label{thm:plummer}
    Every $3$-connected claw-free graph of even order is a brick.
\end{theorem}

% \begin{proof}
%     Suppose that $G$ is bicritical. Let $S \subseteq V(G)$ with $|S| \ge 2$, and choose two distinct vertices $x,y \in S$. Since $G-\{x,y\}$ admits a perfect matching, Tutte's theorem implies that $$o(G-S) = o((G-\{x,y\})-(S\setminus\{x,y\})) \le |S| - 2.$$
%     Conversely, suppose that $o(G-S) \le |S|-2$ for any $S \subseteq V(G)$ with $|S| \ge 2$. Let $x$ and $y$ be any two distinct vertices of $G$. For every subset $T$ of $V(G)\setminus\{x,y\}$, let $S = T\cup\{x,y\}$. Then, $|S| \ge 2$ and $o((G-\{x,y\})-T) = o(G-S) \le |S|-2 = |T|$. Thus $G-\{x,y\}$ satisfies Tutte's condition and hence admits a perfect matching. Therefore $G$ is bicritical.
% \end{proof}

A graph is {\em minimally $3$-connected} if it is $3$-connected but deleting any edge destroys $3$-connectivity.

\begin{theorem}[Halin~\cite{Halin1969}]\label{thm:Halin}
    Every cycle of a minimally $3$-connected graph contains at least two cubic vertices.
\end{theorem}

The following lemma characterizes edges whose deletion preserves $3$-connectivity but destroys bicriticality in claw-free bricks.

\begin{lemma}\label{lem:claw-free-nonbicritical}
    Let $G$ be a claw-free brick and let $ab \in E(G)$. Then, $G-ab$ is $3$-connected but not bicritical if and only if $d_G(a) = d_G(b) = 4$ and $N_G(a) \setminus \{b\} = N_G(b)\setminus\{a\}$.
    When these equivalent conditions hold, set $S:= N_G(a) \setminus \{b\} = N_G(b)\setminus\{a\}$. Then $E(G[S]) \ne \emptyset$, and for every nonisolated vertex $x$ of $G[S]$, both $G-ax$ and $G-bx$ are bricks.
\end{lemma}
\begin{proof}
Let $G$ be a claw-free brick and $ab \in E(G)$, and let $G':= G-ab$.

We first prove the ``if'' part. Suppose that $d_G(a) = d_G(b) = 4$ and $N_G(a) \setminus \{b\} = N_G(b)\setminus\{a\}$. Let $S$ denote this common set. Then, $|S| = d_G(a) -1 = d_G(b)-1 = 3$. We claim that $G'$ is $3$-connected. Otherwise, $G'$ has a vertex cut $X \subseteq V(G')$ with $|X| \le 2$. 
Since $G$ is $3$-connected, $a, b$ must be in different components of $G'-X$, which implies that $S = N_G(a) \cap N_G(b) \subseteq X$ and so $|X| \ge 3$, a contradiction. Thus $G'$ is $3$-connected. 
Since $|G'|=|G|$ is even and $|S| = 3$, it follows that $o(G'-S) \equiv |S| \equiv 1 \pmod 2$. Since both $a, b$ are isolated vertices in $G'-S$, we have $o(G'-S) \ge 2$, and hence $o(G'-S) \ge 3 > |S| - 2$. Thus $G'$ is not bicritical by Lemma~\ref{lem:tutte-bicritical-fact}. 

Conversely, suppose that $G'$ is $3$-connected but not bicritical. By Lemma~\ref{lem:tutte-bicritical-fact}, there is a set $S\subseteq V(G') = V(G)$ with $|S| \ge 2$ satisfying $o(G'-S) \ge |S|-1$. Since $|G'|$ is even, $o(G'-S)\equiv |S|\pmod 2$, and so $o(G'-S)\ge |S|$. On the other hand, adding the edge $ab$ to $G'-S$ can decrease the number of odd components by at most two. Since $G$ is bicritical, by Lemma~\ref{lem:tutte-bicritical-fact}, we have
\[
|S|-2\ge o(G-S) \ge o(G'-S)-2 \ge |S|-2.
\]
Therefore,
\[
o(G-S)=|S|-2 \quad\text{and}\quad o(G'-S)=|S|.
\]
In particular, adding the edge $ab$ to $G'-S$ reduces the number of odd components by exactly $2$. So $a$ and $b$ lie in two distinct odd components of $G'-S$. Denote by $A$ and $B$ the odd components containing $a$ and $b$, respectively. Since $G'$ is $3$-connected and $G'-S$ is disconnected, we have $|S| \ge 3$, and every component of $G'-S$ has at least three distinct neighbors in $S$.

Now count the pairs $(s, C)$, where $s \in S$, $C$ is a component of $G'-S$, and $s$ has a neighbor in $C$. Since $G'-S$ has exactly $|S|$ odd components, and possibly some even components, there are at least $3|S|$ such pairs. On the other hand, since $G$ is claw-free, each vertex of $S$ contributes to at most three such pairs, giving a total of at most $3|S|$. Hence equality holds. Therefore $G'-S$ has no even component, every component of $G'-S$ has exactly three distinct neighbors in $S$, and every vertex of $S$ has neighbors in exactly three components of $G'-S$.  

We claim that, for every $s \in S$, $N_G(s) \cap V(A) = \{a\}$ and $N_G(s) \cap V(B) = \{b\}$. Indeed, choose one neighbor of $s$ from each of the three odd components in which $s$ has neighbors. As $G$ is claw-free, these three vertices cannot be pairwise nonadjacent in $G$. Since the only possible edge of $G$ joining distinct components of $G'-S$ is $ab$, two of the chosen vertices must be $a$ and $b$. Thus $s$ is adjacent to both $a$ and $b$. If $s$ has another neighbor $a' \in V(A)\setminus \{a\}$, then $a'$, $b$, and a neighbor of $s$ in the third component would be pairwise nonadjacent in $G$, inducing a claw centered at $s$. Therefore $N_G(s) \cap V(A) = \{a\}$. Similarly, $N_G(s) \cap V(B) = \{b\}$. This also implies that if $V(A) \setminus \{a\} \ne \emptyset$, then deleting $a$ would disconnect $G'$, contradicting its $3$-connectivity. Hence $V(A) = \{a\}$, and similarly $V(B) = \{b\}$. Since every vertex of $S$ is adjacent to $a$ and $A=\{a\}$ has exactly three distinct neighbors in $S$, it follows that $|S| = 3$. Thus $N_G(a) \setminus \{b\} = N_G(b) \setminus \{a\} = S$ and $d_G(a) = d_G(b) = 4$. This proves the equivalence.

Suppose now that the equivalent conditions above hold. The vertex $a$ is adjacent to all three vertices of $S$, and since $G$ is claw-free, $E(G[S]) \ne \emptyset$. Write $S = \{x,y,z\}$ and assume that $xy \in E(G[S])$. We claim that $G-ax$ is a brick. By the necessity argument above, $G'-S$ has exactly three components, namely $\{a\}$, $\{b\}$, and one odd component, say $C_0$. Then, $V(G) = \{a, b\} \dot\cup S \dot\cup V(C_0)$, and every vertex of $S$ has a neighbor in $C_0$. 

We now show that $G-ax$ is $3$-connected. Suppose, to the contrary, that $G-ax$ is not $3$-connected. Since $G$ is $3$-connected, there exists $X \subseteq V(G)$ with $|X| \le 2$ such that $a$ and $x$ lie in distinct components of $(G-ax)-X$. Then $N_G(a)\cap N_G(x) \subseteq X$, forcing $X = \{b, y\}$. However, $(G-ax)-X$ is connected, since $C_0$ is connected, both $x$ and $z$ have neighbors in $C_0$, and $a$ is adjacent to $z$. Hence $G-ax$ is $3$-connected. Since $x$ has a neighbor in $C_0$ and $a$ does not, $N_G(a) \setminus \{x\} \neq N_G(x) \setminus \{a\}$. By the necessity direction proved above, this implies that $G-ax$ is bicritical, and therefore it is a brick. The same argument applies to $G-bx$. 
\end{proof}

Combining Lemma~\ref{lem:claw-free-nonbicritical} with Theorem~\ref{thm:plummer}, we obtain the following reduction.

\begin{lemma}\label{lem:reduction}
    The claw-free minimal bricks are precisely the claw-free minimally $3$-connected graphs of even order.
\end{lemma}
\begin{proof}
Let $G$ be a claw-free minimal brick. Since $G$ is bicritical, it has even order. We claim that $G$ is minimally $3$-connected.
Otherwise, $G-e$ is $3$-connected for some edge $e\in E(G)$. Since $G$ is a minimal brick, $G-e$ is not a brick and hence is not bicritical. By Lemma~\ref{lem:claw-free-nonbicritical},
there exists an edge $f\in E(G)$ such that $G-f$ is a brick, contradicting the minimality of $G$. Thus $G$ is minimally $3$-connected.

On the other hand, let $G$ be a claw-free minimally $3$-connected graph of even order. By Theorem~\ref{thm:plummer}, $G$ is a brick. Since $G-e$ is not $3$-connected for any $e\in E(G)$, it cannot be a brick. Hence $G$ is a minimal brick.
\end{proof}

We now describe the structural properties of claw-free minimal bricks. For a graph $G$, let $$V_k := \{v \in V(G) : d_G(v) = k\}.$$

\begin{lemma}\label{lem:structures-3}
Let $G$ be a claw-free minimal brick other than $W_5$. Then, every vertex of $G$ has degree either $3$ or $4$, and the vertices of degree $4$ form an independent set. Moreover, $G[N_G(u)] \cong 2K_2$ for every vertex $u$ of degree $4$, and every cubic vertex is adjacent to at most one vertex of degree $4$. 
\end{lemma}
\begin{proof}
    Let $V_{\ge 4}:=V(G)\setminus V_3$.
    By Lemma~\ref{lem:reduction}, $G$ is minimally $3$-connected. Thus, by Theorem~\ref{thm:Halin}, every cycle of $G$ contains at least two cubic vertices, and hence $G[V_{\ge 4}]$ is a forest. In particular, every triangle of $G$ contains at least two cubic vertices. Since $G$ is claw-free, every cubic vertex of $G$ belongs to a triangle, which implies that every cubic vertex has a cubic neighbor.
    
    We first show that $V_{\ge 4}$ is independent. Suppose for a contradiction that $G[V_{\ge 4}]$ contains an edge. Choose a leaf $u \in V_{\ge 4}$ in a nontrivial component of $G[V_{\ge 4}]$, and let $w$ be its unique neighbor in $V_{\ge 4}$. As $d_G(u) \ge 4$, $u$ has at least three neighbors in $V_3$. Since every triangle of $G$ contains at least two cubic vertices, no cubic neighbor of $u$ is adjacent to $w$. Thus any two cubic neighbors of $u$ must be adjacent; otherwise, they, together with $w$, would form an induced claw centered at $u$, a contradiction. Hence the cubic neighbors of $u$ form a clique. Since each of them is cubic and is adjacent to $u$, there can be at most three such vertices. It follows that $|N_G(u) \cap V_3| = 3$. These three vertices form a triangle and have no neighbors outside it other than $u$. Thus deleting $u$ separates this triangle from the rest of the graph, contradicting the $3$-connectivity of $G$. Therefore $V_{\ge 4}$ is independent.

    We claim that $d_G(u) = 4$ and $G[N_G(u)] \cong 2K_2$ for any $u \in V_{\ge 4}$. Since $V_{\ge 4}$ is independent, $N_G(u) \subseteq V_3$. Every vertex of $G[N_G(u)]$ is cubic in $G$ and adjacent to $u$, so $\Delta(G[N_G(u)]) \le 2$. Since $G$ is claw-free, $\alpha(G[N_G(u)]) \le 2$. Observe that $G[N_G(u)]$ is triangle-free. Indeed, if three vertices of $N_G(u)$ form a triangle, then deleting $u$ would separate this triangle from another neighbor of $u$, contradicting the $3$-connectivity of $G$. Thus every component of $G[N_G(u)]$ is a path or a cycle of length at least $4$. If $d_G(u) \ge 6$, then $\alpha(G[N_G(u)]) \ge 3$, a contradiction. If $d_G(u) = 5$, then the only possibility is $G[N_G(u)] \cong C_5$. Since the vertices of $C_5$ are cubic in $G$, this forces $G \cong W_5$, a contradiction. Therefore $d_G(u) = 4$, and $G[N_G(u)]$ is isomorphic to one of $2K_2$, $P_4$, and $C_4$. If $G[N_G(u)] \cong C_4$, then $G \cong W_4$, which has odd order, contradicting the bicriticality of $G$. If $G[N_G(u)] \cong P_4$, deleting the two endvertices of the path separates $u$ and the two internal vertices from the rest of the graph, which is nonempty because $G$ has even order. This contradicts the $3$-connectivity of $G$. Hence $V_{\ge 4} = V_4$ and $G[N_G(u)] \cong 2K_2$ for every $u \in V_4$.

    We now show that every vertex of $V_3$ has at most one neighbor in $V_4$. Suppose, to the contrary, that some vertex $v \in V_3$ satisfies $N_G(v) = \{u, w, z\}$, where $u, w \in V_4$. Then $uw \notin E(G)$. Since every cubic vertex has a cubic neighbor, we have $z \in V_3$. As $G[N_G(u)] \cong 2K_2$, the unique neighbor of $v$ in $N_G(u)$ must be $z$, so $uz \in E(G)$. Similarly, $wz \in E(G)$. Since $z$ is cubic, $N_G(z) = \{u,v,w\}$. Thus deleting $u$ and $w$ separates $\{v,z\}$ from the rest of the graph, which is nonempty since $d_G(u) = d_G(w) = 4$, contradicting the $3$-connectivity of $G$.
\end{proof}

\begin{lemma}\label{lem:triangle}
Let $G$ be a claw-free minimal brick other than $K_4$ and $W_5$. The following hold: 
\begin{itemize}
    \item [{\em (a)}] Every cubic vertex lies in exactly one triangle;
    \item [{\em (b)}] Each noncubic vertex lies in exactly two triangles, both containing two cubic vertices.
\end{itemize}     
\end{lemma}
\begin{proof}
    Part $(b)$ is immediate from Lemma~\ref{lem:structures-3}, so it remains to prove $(a)$.
    Since $G$ is claw-free, every cubic vertex of $G$ lies in a triangle. Suppose, to the contrary, that $u \in V_3$ with $N_G(u) = \{v,x,y\}$ lies in two triangles. There are at least two edges in $G[N_G(u)]$ and they share a common vertex. We, without loss of generality, assume that $xv, yv \in E(G)$. Since $u$ has degree $2$ in $G[N_G(v)]$, Lemma~\ref{lem:structures-3} implies that $d_G(v) = 3$, and hence $N_G(v) = \{u,x,y\}$. Thus deleting $x$ and $y$ separates $\{u,v\}$ from the rest of $G$, unless $V(G) = \{u,v,x,y\}$. The former contradicts the $3$-connectivity of $G$, while the latter forces $G\cong K_4$, again a contradiction. Therefore every cubic vertex lies in exactly one triangle. 
\end{proof}

\section{Characterization of claw-free minimal bricks}\label{sec:pf-main}

In this section, we prove the characterization of claw-free minimal bricks and determine their possible orders. For a loopless multigraph $Q$ and $M\subseteq E(Q)$, let $R(Q,M)$ be obtained from $Q$ by subdividing every edge in $E(Q)\setminus M$ exactly once and leaving the edges of $M$ unsubdivided. Define
\[
    G(Q,M):=L(R(Q,M)).
\]

\begin{proof}[Proof of Theorem~\ref{thm:main1-char}]
We first show that every graph described in the theorem is a claw-free minimal brick. Clearly, $K_4$ and $W_5$ are claw-free minimal bricks. 
Let $Q$ be a loopless $3$-edge-connected cubic multigraph, and let $M \subseteq E(Q)$ be a matching of even size. Set $R:= R(Q, M)$ and $G:= G(Q, M)$.

We observe that $Q$ is simple unless $|Q| = 2$. Indeed, if $|Q| \ge 4$ and two vertices $u, v$ are joined by parallel edges, then the edge boundary of $\{u,v\}$ has size at most two, contradicting the $3$-edge-connectivity of $Q$. Since $Q$ is cubic, $|Q|$ is even. Hence if $Q$ is not simple, then $|Q| = 2$, in which case $Q$ consists of two vertices joined by three parallel edges. As $M$ is a matching of even size, $M = \emptyset$. Consequently, $R \cong K_{2,3}$ and $G$ is the triangular prism. In either case, $R$ is simple, and hence so is $G = L(R)$. Since every line graph is claw-free, $G$ is claw-free. Moreover, $|G| = |E(R)| = 2|E(Q)| - |M|$, which is even because $|M|$ is even. Therefore $G$ is a claw-free simple graph of even order.

We claim that $G$ is $3$-connected. Let $X \subseteq V(G)$ with $|X| \le 2$, and let $F \subseteq E(R)$ be the corresponding set of edges. Then, $|F| = |X| \le 2$ and $G-X = L(R-F)$. The edges of $F$ lie on the paths of $R$ corresponding to at most two edges of $Q$. Since $Q$ is $3$-edge-connected, deleting these edges leaves $Q$ connected. Hence all original vertices of $Q$ lie in one component of $R-F$. Any subdivision vertex outside this component must be isolated, since any surviving incident edge would join it to an original vertex of $Q$. Thus $R-F$ has at most one component containing edges. Therefore $L(R-F) = G-X$ is connected. It follows that $G$ is $3$-connected. 

We now show that $G$ is minimally $3$-connected. By construction, each edge of $M$ corresponds to a vertex of degree $4$ in $G$, while the two edges arising from each subdivided edge of $E(Q)\setminus M$ correspond to cubic vertices of $G$. Since $M$ is a matching, no two edges of $M$ share a common vertex in $Q$, and hence in $R$. Thus no two vertices of degree $4$ in $G$ are adjacent. This implies that every edge of $G$ is incident with a cubic vertex. Consequently, deleting any edge of $G$ leaves a vertex of degree $2$, so the resulting graph is not $3$-connected. Thus $G$ is minimally $3$-connected. Therefore $G$ is a claw-free minimally $3$-connected graph of even order, and hence is a claw-free minimal brick by Lemma~\ref{lem:reduction}.

Conversely, let $G$ be a claw-free minimal brick other than $K_4$ and $W_5$. By Lemma~\ref{lem:structures-3}, every vertex of $G$ has degree $3$ or $4$. If $G$ is cubic, set $T:=G$. Otherwise, choose a vertex $v$ of degree $4$. Since $G[N_G(v)] \cong 2K_2$, write $N_G(v) = \{u_1,u_2, w_1, w_2\}$, where $u_1u_2, w_1w_2 \in E(G)$.  
Split $v$ into two adjacent vertices $v_1$ and $v_2$ such that $v_1$ is adjacent to $u_1, u_2$, while $v_2$ is adjacent to $w_1, w_2$. Denote the resulting graph by $G_v$. Both $v_1$ and $v_2$ are cubic, and the degrees of all other vertices remain unchanged.
We claim that $G_v$ is $3$-connected. Let $X \subseteq V(G_v)$ with $|X| \le 2$. If $X \cap \{v_1, v_2\} = \emptyset$, then contracting the edge $v_1v_2$ in $G_v-X$ yields $G-X$. Since $G$ is $3$-connected, $G-X$ is connected, so is $G_v-X$. If $X = \{v_1,v_2\}$, then $G_v-X = G-v$, which is also connected. It remains to consider the case where $X$ contains exactly one of $\{v_1, v_2\}$. By symmetry, assume that $v_1\in X$. Set $Y:= X\setminus\{v_1\}$, so $|Y| \le 1$. Clearly, $G-(Y\cup\{v\})= (G_v-X)-v_2$. Since $G-(Y\cup\{v\})$ is connected and $v_2$ is adjacent to at least one vertex of $G-(Y\cup\{v\})$, the graph $G_v-X$ is connected. Thus $G_v$ is $3$-connected. Since the vertices of degree $4$ in $G$ are pairwise nonadjacent, splitting $v$ does not change the neighborhood of any other vertex of degree $4$. We repeat this splitting operation at all remaining vertices of degree $4$, and denote the resulting graph by $T$. Each operation preserves $3$-connectivity by the same argument. Therefore, in either case, $T$ is simple, cubic and $3$-connected.

We then show that the vertices of $T$ can be partitioned into triangles. Every vertex of $T$ is either a cubic vertex of $G$ or a vertex created by splitting a vertex of degree $4$ in $G$. By Lemma~\ref{lem:triangle}$(a)$, every cubic vertex of $G$ lies in exactly one triangle. By Lemma~\ref{lem:triangle}$(b)$, every vertex of degree $4$ in $G$ lies in exactly two triangles, each containing two cubic vertices. Splitting such a vertex preserves its two triangles, with the original vertex replaced by a distinct new vertex in each triangle. Thus every vertex of $T$ lies in exactly one triangle. Let $\mathcal{T}$ denote the collection of these triangles. 

We now construct a loopless $3$-edge-connected cubic multigraph $Q$ and a matching $M\subseteq E(Q)$ of even size such that $G\cong G(Q,M)$. Contract each triangle in $\mathcal{T}$ into a single vertex, retaining all edges joining distinct triangles, and denote the resulting multigraph by $Q$. Since $T$ is cubic and its vertices are partitioned into triangles, $Q$ is loopless and cubic. Moreover, $Q$ is $3$-edge-connected, since otherwise an edge cut of $Q$ of size at most $2$ would correspond to an edge cut of the same size in $T$, contradicting the $3$-connectivity of $T$. Let $M$ be the set of edges of $T$ introduced by the splitting operations, which is empty if $G$ is cubic. Since these edges join distinct triangles of $T$, they remain as edges of $Q$. Thus $M \subseteq E(Q)$. Since each triangle of $G$ contains at most one vertex of degree $4$, each triangle of $T$ contains at most one newly created vertex. Consequently, no two edges of $M$ share a common vertex in $Q$, so $M$ is a matching. To show that $|M|$ is even, observe that $3|V_3|+4|V_4|=2|E(G)|$, so $|V_3|$ is even. Since $G$ is bicritical, $|G| = |V_3|+|V_4|$ is also even. Hence $|M| = |V_4|$ is even.

It remains to show that $G\cong G(Q,M)$. For each vertex of $Q$, the three incident edges of $R(Q,M)$ form a triangle in $L(R(Q,M))$, corresponding to a triangle of $T$. Each edge in $E(Q)\setminus M$ is subdivided once in $R(Q,M)$, producing two 
adjacent vertices in $L(R(Q,M))$, one in each of the corresponding triangles. This recovers the edges of $T$ joining distinct triangles outside $M$. On the other hand, each edge of $M$ remains unsubdivided in $R(Q,M)$, so the corresponding triangles share a vertex in $L(R(Q,M))$. This identifies the two vertices of $T$ created by the corresponding splitting operation. Thus $L(R(Q,M))$ is obtained from $T$ by contracting the edges introduced by the splitting operations, and hence $G(Q,M)=L(R(Q,M))\cong G$.
\end{proof}

\begin{corollary}
    The orders of claw-free minimal bricks are precisely
    $$
    \{4,6,10,12\} \cup \{n \ge 16 : n \text{ is even}\}.
    $$
    In particular, there is no claw-free  minimal brick of order $8$ or $14$.
\end{corollary}
\begin{proof}
    Let $G$ be a claw-free minimal brick of order $n$ other than $K_4$ and $W_5$. By Theorem~\ref{thm:main1-char}, there exist a loopless cubic $3$-edge-connected multigraph $Q$ and a matching $M \subseteq E(Q)$ of even size such that $G \cong G(Q,M)$. Let $|Q| = 2k$ and $|M| = s$. Then, $s$ is even, $0 \le s \le k$, and $$n = |E(R(Q, M))| = 2|E(Q)|-|M| = 6k-s.$$
    If $k=1$, then $s = 0$ and $n = 6$. 
    If $k = 2$, then $s \in \{0, 2\}$, giving $n \in \{10, 12\}$. 
    If $k \ge 3$, then $5k \le n \le 6k$. Since $n$ is even, we have $n \ge 16$. Together with the exceptional cases $K_4$ and $W_5$, this shows that the order of every claw-free minimal brick belongs to $\{4, 6, 10, 12\} \cup \{n \ge 16 : n \text{ is even}\}$.

    Conversely, let $n \in \{6,10,12\} \cup \{n \ge 16  : n \text{ is even}\}$. Set $k = \lceil n/6 \rceil$ and $s = 6k-n$. Then $s$ is even and $0 \le s \le k$. We show that there exists a loopless $3$-edge-connected cubic multigraph $Q$ of order $2k$ with a matching $M$ of size $s$.

    If $k = 1$, take $Q$ to be the multigraph consisting of two vertices joined by three parallel edges and let $M = \emptyset$. 
    If $k = 2$, take $Q = K_4$ and let $M$ be a matching of size $s$. If $k \ge 3$, take the \emph{$k$-prism graph} $Q = C_k \square K_2$, which consists of two disjoint $k$-cycles joined by a perfect matching of $k$ edges. The edges of this matching are called the \emph{rungs}. Let $M$ consist of any $s$ rungs.
    
    In each case, Theorem~\ref{thm:main1-char} implies that $G(Q,M)$ is a claw-free minimal brick of order $6k-s = n$. Together with $K_4$, this completes the proof.
\end{proof}

\section{Applications}\label{sec:applications}
In this section, we show Theorems~\ref{thm:counting} and~\ref{thm:edge-classes}. Let $G$ be a claw-free minimal brick of order $n$ other than $K_4$ and $W_5$. By Theorem~\ref{thm:main1-char}, there exist a loopless  $3$-edge-connected cubic multigraph $Q$ and a matching $M$ of even size such that $G \cong G(Q,M)$. Let $|Q| = 2k$ and $|M| = s$. Then, $s$ is even and $0 \le s \le k$. 
Each edge of $M$ corresponds to a vertex of degree $4$ in $G$, while each edge of $E(Q) \setminus M$ gives rise to two cubic vertices of $G$. It follows that 
$$
n = 6k-s \quad \text{and} \quad |V_3| = 6k-2s.
$$
For each $v \in V(Q)$, let $T_v$ be the triangle of $G$ corresponding to the three edges of $R(Q,M)$ incident with $v$. These triangles are pairwise edge-disjoint, since $R(Q,M)$ is simple. Every triangle of $G$ contains a cubic vertex, and every cubic vertex belongs to some $T_v$. Hence by Lemma~\ref{lem:triangle}$(a)$, the triangles $T_v$, $v\in V(Q)$, are precisely the triangles of $G$.

\subsection{Cubic vertices}\label{subsec:cubic-counting}

\begin{proof}[Proof of Theorem~\ref{thm:counting}]
For $G = K_4$ and $G = W_5$, the subgraph $G[V_3]$ is $K_4$ and $C_5$, respectively, and the result follows directly. Now assume that $G$ is a claw-free minimal brick other than $K_4$ and $W_5$, and use the notation above. 

We claim that $G[V_3]$ admits a perfect matching $P$ such that $G[V_3]-P$ is a disjoint union of triangles and single edges. For each edge $e\in E(Q)\setminus M$, the subdivision of $e$ in $R(Q,M)$ produces two adjacent edges, which correspond to two adjacent cubic vertices of $G$. Let $P$ be the set of all edges of $G$ joining these pairs. Since every cubic vertex of $G$ belongs to exactly one such pair, $P$ is a perfect matching of $G[V_3]$. We now determine the components of $G[V_3]-P$. By construction, every edge of $G$ outside $P$ belongs to exactly one triangle $T_v$. If $v$ is not incident with an edge of $M$, then all three vertices of $T_v$ are cubic, so $T_v$ is a triangle in $G[V_3]-P$. Otherwise, exactly two vertices of $T_v$ are cubic, and the edge joining them is in $G[V_3]-P$. Since every cubic vertex of $G$ belongs to exactly one triangle, these triangles and single edges are pairwise vertex-disjoint, and hence form the components of $G[V_3]-P$. Since $M$ is a matching of size $s$, exactly $2s$ vertices of $Q$ are incident with edges of $M$. Thus $G[V_3]-P$ consists of $2k-2s$ triangles and $2s$ single edges. 

Since $|V_3| = 6k-2s$, we have $$|P| = 3k-s \quad \text{and} \quad |E(G[V_3])| = 3(2k-2s)+2s+(3k-s) = 9k-5s.$$ 
As $s \le k$ and $n = 6k-s$, we obtain $$|E(G[V_3])| - |V_3| = 3(k-s) \ge 0 \quad \text{and} \quad 5|V_3| - 4n = 6(k-s) \ge 0.$$ Therefore $|E(G[V_3])| \ge |V_3| \ge 4n/5$.

We now show that both inequalities are sharp. Let $k \ge 4$ be even, let $Q=C_k\square K_2$, and let $M$ consist of all $k$ rungs of $Q$. By Theorem~\ref{thm:main1-char}, the graph $G := G(Q,M)$ is a claw-free minimal brick, and since $s = |M| = k$, it follows that $|E(G[V_3])|=|V_3|=4|G|/5$.
\end{proof}

\subsection{Removable, $b$-invariant and thin edges}
\label{subsec:edge-classes}
Recall that $R(G)$, $I(G)$, and $T(G)$ denote the sets of removable, $b$-invariant, and thin edges of a graph $G$, respectively. We also recall the following result of Xue and Lu.

\begin{theorem}[Xue and Lu~\cite{XL2025thin}]\label{thm:XL}
Let $G$ be a claw-free brick of order at least $8$. Then every $b$-invariant edge of $G$ is thin, and every removable edge incident with a cubic vertex is $b$-invariant.
\end{theorem}

\begin{proof}[Proof of Theorem~\ref{thm:edge-classes}] 
Let $G$ be a claw-free minimal brick of order at least $8$, and use the notation introduced at the beginning of this section. By Theorem~\ref{thm:main1-char}, $Q$ is simple, and by Lemma~\ref{lem:structures-3}, every edge of $G$ is incident with a cubic vertex, so $R(G) = I(G) = T(G)$ by Theorem~\ref{thm:XL}. We now show that the removable edges of $G$ are precisely those lying in a triangle.

Let $P$ be the perfect matching constructed in the proof of Theorem~\ref{thm:counting}. Recall that every edge of $G$ belongs either to $P$ or to a unique triangle $T_v$ for some $v \in V(Q)$. We first show that every edge of $P$ is nonremovable. Let $v_1v_2 \in P$ and $T_v$ be the triangle containing $v_1$. Write $V(T_v) = \{v_1, a, b\}$. Since $v_1$ is cubic, every perfect matching of $G-v_1v_2$ contains either $v_1a$ or $v_1b$. Thus $ab$ belongs to no perfect matching of $G-v_1v_2$, so $v_1v_2$ is nonremovable. 

We then show that every edge of $T_v$ is removable for each $v \in V(Q)$. We observe that, for every $w \in V(Q)$ and $e\in E(Q)$, both $Q-w$ and $Q-e$ are connected and bridgeless. Indeed, suppose that $Q-w$ is disconnected or has a bridge $f$. In the former case, each component of $Q-w$ is joined to $w$ by at least three edges, while in the latter case, each component of $(Q-w)-f$ is joined to $w$ by at least two edges, since $Q$ is $3$-edge-connected. In either case, this contradicts $d_Q(w) = 3$. For $Q-e$, this conclusion follows directly from the $3$-edge-connectivity of $Q$. It follows that every vertex of $Q-w$ and $Q-e$ lies on a cycle. 

We use the following correspondence between perfect matchings of $G$ and orientations of $R(Q,M)$ in which every vertex has even indegree. Since $G \cong L(R(Q,M))$, each perfect matching of $G$ partitions $E(R(Q,M))$ into pairs of incident edges. Orient two edges in each pair toward their common vertex. Since $\Delta(R(Q,M)) \le 3$, every vertex has indegree $0$ or $2$. Conversely, any orientation of $R(Q,M)$ in which every vertex has indegree $0$ or $2$ determines a perfect matching of $G$ by pairing two incoming edges at each vertex of indegree $2$, since these pairs partition $E(R(Q,M))$. Given such an orientation, reversing the direction of all edges on a cycle of $R(Q,M)$ preserves the parity of the indegree at every vertex, and hence yields another perfect matching of $G$.

Now fix $v \in V(Q)$ and $e \in E(T_v)$. To show that $e$ is removable, it suffices to show that every edge of $G-e$ belongs to a perfect matching of $G-e$, since $G-e$ is connected. Let $f \in E(G-e)$. Since $G$ is matching covered, there is a perfect matching $F$ of $G$ containing $f$. If $e \notin F$, then $F$ is also a perfect matching of $G-e$ containing $f$, and we are done. We may assume that $e \in F$. Then $f \notin E(T_v)$, since two edges of a triangle cannot belong to a matching. If $f \in E(T_w)$ for some $w \ne v$, choose a cycle through $v$ in $Q-w$. 
Otherwise, $f \in P$ corresponds to an edge $f_0\in E(Q)\setminus M$, and we choose a cycle through $v$ in $Q-f_0$.
In either case, reversing the corresponding cycle in $R(Q,M)$ in the orientation determined by $F$ yields another perfect matching. Since the cycle passes through $v$, the resulting perfect matching avoids $e$. Moreover, since the cycle avoids $w$ or $f_0$, it still contains $f$. Thus every edge of $G-e$ belongs to a perfect matching of $G-e$, and hence $e$ is removable.

There are $2k$ pairwise edge-disjoint triangles $T_v$ in $G$, $|G|=6k-s$, and $|V_3| = 6k-2s$. Thus $$|R(G)|=|I(G)|=|T(G)|=6k=2|G|-|V_3|.$$
Since $|V_3| \le |G|$ and $|V_3| \ge 4|G|/5$ by Theorem~\ref{thm:counting}, the bounds follow. The lower bound is attained when $G$ is cubic, while the upper bound is attained by the construction in the proof of Theorem~\ref{thm:counting}.
\end{proof}

\section*{Conflict of Interest Statement}
The authors declare that they have no conflict of interest regarding the publication of this paper.

\bibliographystyle{plain}
\bibliography{reference}

@article{NORINE2006,
title = {Minimal bricks},
journal = {Journal of Combinatorial Theory, Series B},
volume = {96},
number = {4},
pages = {505-513},
year = {2006},
issn = {0095-8956},
doi = {https://doi.org/10.1016/j.jctb.2005.10.005},
url = {https://www.sciencedirect.com/science/article/pii/S0095895605001516},
author = {S. Norine and R. Thomas}
}

@article{LZL2014,
author = {Lin, F. and Zhang, L. and Lu, F.},
title = {The Cubic Vertices of Minimal Bricks},
journal = {Journal of Graph Theory},
volume = {76},
number = {1},
pages = {20-33},
doi = {https://doi.org/10.1002/jgt.21747},
url = {https://onlinelibrary.wiley.com/doi/abs/10.1002/jgt.21747},
eprint = {https://onlinelibrary.wiley.com/doi/pdf/10.1002/jgt.21747},
year = {2014}
}

@article{DECARVALHO2006,
title = {How to build a brick},
journal = {Discrete Mathematics},
volume = {306},
number = {19},
pages = {2383-2410},
year = {2006},
note = {Creation and Recreation: A Tribute to the Memory of Claude Berge},
issn = {0012-365X},
doi = {https://doi.org/10.1016/j.disc.2005.12.032},
url = {https://www.sciencedirect.com/science/article/pii/S0012365X06003281},
author = {M. H. {de Carvalho} and C. L. Lucchesi and U.S.R. Murty}
}

@article{Feng2024,
  title={On a Norine-Thomas conjecture concerning minimal bricks},
  author={X. Feng},
  journal={Journal of Graph Theory},
  year={2024},
  volume={108},
  pages={162-172},
  url={https://api.semanticscholar.org/CorpusID:272636016}
}

@article{HE2024,
title = {The cubic vertices of solid minimal bricks},
journal = {Discrete Mathematics},
volume = {347},
number = {2},
pages = {113746},
year = {2024},
issn = {0012-365X},
doi = {https://doi.org/10.1016/j.disc.2023.113746},
url = {https://www.sciencedirect.com/science/article/pii/S0012365X23004326},
author = {X. He and F. Lu},
}

@article{Plummer1994,
  title={Extending matchings in claw-free graphs},
  author={Plummer, M. D.},
  journal={Discrete mathematics},
  volume={125},
  number={1-3},
  pages={301--307},
  year={1994},
  publisher={Elsevier}
}

@article{Halin1969,
  author  = {Halin, R.},
  title   = {Untersuchungen {\"u}ber minimale $n$-fach zusammenh{\"a}ngende Graphen},
  journal = {Mathematische Annalen},
  volume  = {182},
  year    = {1969},
  pages   = {175-188},
  doi     = {10.1007/BF01350321}
}

@article{XL2025thin, 
title={Thin Edges in Claw-Free Bricks}, 
volume={32}, 
url={https://www.combinatorics.org/ojs/index.php/eljc/article/view/v32i3p17}, 
DOI={10.37236/12976},
number={3}, 
journal={The Electronic Journal of Combinatorics}, 
author={Xue, J. and Lu, F.}, 
year={2025}, 
pages={3-17} }

@inproceedings{Yuliu2009,
  title={Graph Factors and Matching Extensions},
  author={Q. Yu and G. Liu},
  year={2009},
  url={https://api.semanticscholar.org/CorpusID:117031107}
}

@article{WLZ2024,
  title={Removable edges in claw-free bricks},
  author={Wu, X. and Lu, F. and Zhang, L.},
  journal={Graphs and Combinatorics},
  volume={40},
  number={3},
  pages={43},
  year={2024},
  publisher={Springer}
}

@article{ZHANG2025,
title = {Claw-free minimal matching covered graphs},
journal = {Discrete Applied Mathematics},
volume = {370},
pages = {11-21},
year = {2025},
issn = {0166-218X},
doi = {https://doi.org/10.1016/j.dam.2025.02.045},
url = {https://www.sciencedirect.com/science/article/pii/S0166218X25001155},
author = {Y. Zhang and X. Wang and J. Yuan and C. T. Ng and T. C. E. Cheng}
}

@book{bondy2011,
  title={Graph Theory},
  author={Bondy, A. and Murty, U. S. R.},
  isbn={9781846289699},
  lccn={2007923502},
  series={Graduate Texts in Mathematics},
  url={https://books.google.com/books?id=HuDFMwZOwcsC},
  year={2011},
  publisher={Springer London}
}

@article{CFL2026,
author = {Chen, G. and Feng, X. and Lu, F.},
title = {Existence of Adjacent Cubic Vertices in Minimal Bricks},
journal = {Journal of Graph Theory},
pages = {},
doi = {10.1002/jgt.70134},
url = {https://onlinelibrary.wiley.com/doi/abs/10.1002/jgt.70134},
eprint = {https://onlinelibrary.wiley.com/doi/pdf/10.1002/jgt.70134},
note = {Published online at \url{https://doi.org/10.1002/jgt.70134}}
}

@article{LOVASZ1987,
title = {Matching structure and the matching lattice},
journal = {Journal of Combinatorial Theory, Series B},
volume = {43},
number = {2},
pages = {187-222},
year = {1987},
issn = {0095-8956},
doi = {https://doi.org/10.1016/0095-8956(87)90021-9},
url = {https://www.sciencedirect.com/science/article/pii/0095895687900219},
author = {L. Lov\'asz}
}

@article{CLM2002I,
title = {On a Conjecture of Lov\'asz Concerning Bricks: I. The Characteristic of a Matching Covered Graph},
journal = {Journal of Combinatorial Theory, Series B},
volume = {85},
number = {1},
pages = {94-136},
year = {2002},
issn = {0095-8956},
doi = {https://doi.org/10.1006/jctb.2001.2091},
url = {https://www.sciencedirect.com/science/article/pii/S009589560192091X},
author = {M. H. {de Carvalho} and C. L. Lucchesi and U.S.R. Murty},
}

@article{CLM2002II,
title = {On a Conjecture of Lov\'asz Concerning Bricks: {II}. Bricks of Finite Characteristic},
journal = {Journal of Combinatorial Theory, Series B},
volume = {85},
number = {1},
pages = {137-180},
year = {2002},
issn = {0095-8956},
doi = {https://doi.org/10.1006/jctb.2001.2092},
url = {https://www.sciencedirect.com/science/article/pii/S0095895601920921},
author = {M. H. {de Carvalho} and C. L. Lucchesi and U.S.R. Murty},
}

@article{Edmonds1982,
  title={Brick decompositions and the matching rank of graphs},
  author={Edmonds, J. and Pulleyblank, W. R. and Lov{\'a}sz, L.},
  journal={Combinatorica},
  volume={2},
  number={3},
  pages={247--274},
  year={1982},
  publisher={Springer}
}
\end{document}